\documentclass[10pt,a4paper,reqno]{amsart}
\usepackage{amsfonts,amsthm,latexsym,amsmath,amssymb,amscd,amsmath,
epsf}

\theoremstyle{plain}
\newtheorem{theorem}           {Theorem}

\newtheorem{corollary}         {Corollary}

\theoremstyle{definition}

\theoremstyle{remark}
\newtheorem{remark}            {Remark}

\theoremstyle{plain}

\theoremstyle{definition}

\theoremstyle{remark}

\newtheorem{e-proposition}[theorem]{Proposition}
\newtheorem{e-definition}[theorem]{Definition\rm}
\input cyracc.def
\newfam\cyrfam
\font\tencyr=wncyr10
\font\sevencyr=wncyr7

\textfont\cyrfam=\tencyr \scriptfont\cyrfam=\sevencyr
    \scriptscriptfont\cyrfam=\sevencyr

\font\tencyr=wncyr10

\def\newop#1{\expandafter\def\csname #1\endcsname{\mathop{\rm
#1}\nolimits}}

\newop{slc}

\begin{document}
        \numberwithin{equation}{section}
       
  % \leftline{\eightit 'This result is not just proved, it is also true' -- mathematical folklore }

        \title[On  Schur-Szeg\"o composition]
        {On   Schur-Szeg\"o composition of polynomials}

       \author[V.~Kostov]{Vladimir Kostov}
\address{Laboratoire J.-A.~Dieudonn\'e, UMR 6621 du CNRS, Universit\'e de Nice, Parc Valrose, F-06108 Nice, cedex 2}
\email{kostov@math.unice.fr}

\author[B.~Shapiro]{Boris Shapiro}
       \address{Mathematics Dept., Stockholm University, S-106 91
Stockholm,
       Sweden}
       \email{shapiro@math.su.se}
\date{\today}
\keywords {Schur-Szeg\"o composition, hyperbolic polynomials}
\subjclass{12D10}

\begin{abstract} Schur-Szeg\"o composition of two polynomials of  degree less or equal than a given positive integer  $n$ introduces an interesting semigroup structure on  polynomial spaces  and  is one of the basic  tools in the analytic theory of polynomials, see \cite {RS}.  In the present paper we add several (apparently) new  aspects to the previously known properties of this operation. Namely, we show how it interacts with the stratification of polynomials according to the multiplicities of their zeros and present the induced  semigroup structure on the set of all  ordered partitions of $n$. 
\end{abstract}

\maketitle

%\tableofcontents

The Schur-Szeg\"o composition of  two polynomials 
$P(x)=\sum_{i=0}^nC_n^i{a_ix^i}$ and 
$Q(x)=\sum_{i=0}^nC_n^i{b_ix^i}$ 
is given by 
 $ P\ast Q(x)=\sum_{i=0}^n C_n^i 
{a_ib_i x^i},$   see e.g. \cite {Sz}.
Let  $Pol_n$ denote the linear space of all polynomials in $x$ of degree at 
most $n$. In what follows we always use its  standard monomial basis 
$\mathcal{B}:=(x^n,x^{n-1},\ldots , 1)$. To any polynomial $P\in Pol_n$ 
one can associate the operator $T_P$ which acts diagonally in $\mathcal{B}$ 
and is uniquely determined by the  condition: 
$T_P(1+x)^n=P(x).$ 
Obviously, for  
$P(x)=C_n^0a_0+C_n^1a_1x+\cdots +C_n^na_nx^n$ 
one has  
%\begin{equation}\label{compseq}
$T_P(x^i)={a_i}$, $i=0,1,\ldots ,n$.
%\end{equation} 
Given $P$ as above we refer to  the sequence $\{ a_i\}$ as to 
the {\it diagonal sequence} of $P$.    Any  two such operators $T_P$ and 
$T_Q$ commute and their product  $T_PT_Q$ corresponds in the above sense 
exactly to the Schur-Szeg\"o composition $P\ast Q$. 
The famous composition theorem of Schur and Szeg\"o (see  original \cite {Sz} 
and e.g.  \S 3.4 of \cite {RS} or \S 2 of \cite {CC1}) reads:
\vspace{1mm}

\begin{theorem}\label{th:SS} Given any linear-fractional image 
%ircular domain 
${\mathcal K}$ of the unit disk containing all the roots 
of $P$ one has that any root 
of $P\ast Q$ is the product of some root of $Q$ by $-\gamma$ where 
$\gamma \in \mathcal K$. 
\end{theorem} 
\vspace{1mm}

Geometric consequences of Theorem~\ref{th:SS}, in particular, 
Proposition~\ref{pr:hyper}, can be found in \S~5.5 of \cite {RS}.  
A polynomial $P\in Pol_n$ is 
called {\it hyperbolic} if all its roots are real. Denote by 
$Hyp_n\subset Pol_n$ the set of all hyperbolic polynomials  
and by $Hyp_n^+\subset Hyp_n$ (resp. $Hyp_n^-\subset Hyp_n$)  the set of all 
hyperbolic polynomials with all positive (resp. all negative) roots. 
Denote by $H_{u,v,w}\subset Hyp_n$ 
(where $u,v,w\in {\bf N}\cup 0$, $u+v+w=n$) 
the set of all hyperbolic polynomials 
with $u$ negative and $w$ positive roots and a $v$-fold zero root.
\vspace{1mm}  

\begin{e-proposition}[Theorem 5.5.5 and Corollary 5.5.10 of \cite {RS}]
\label{pr:hyper}
 If  $P, Q\in Hyp _n$ and if $Q\in Hyp _n^+$ or $Q\in Hyp _n^-$, then   
$P\ast Q\in Hyp _n$.  Moreover, all roots of $P\ast Q$ lie in  
$[-M, -m]$ where $M$ is the maximal and $m$ is the minimal 
pairwise product of roots of $P$ and $Q$. 
\end{e-proposition} 
\vspace{1mm}  

%Notice that  an important  class of diagonal sequences was introduced by 
%Polya and Schur \cite {PS}  and later studied in numerous papers.    
A diagonal sequence, (or an operator $T:Pol_n\to Pol_n$ 
acting diagonally in $\mathcal{B}$) is called a {\it finite multiplier 
sequence (FMS)}, see \cite{PS}, 
if it sends $Hyp_n$ into $Hyp_n$.   
The set $\mathcal{M}_n$ of all FMS is a  semigroup. 
For the following characterization of FMS see \cite {CC2}, Theorem 3.7 or 
\cite {CC1}, Theorem 3.1.  
\vspace{1mm}

\begin{theorem}\label{th:MS}
For $T\in$End$(Pol_k^{\mathbb{R}})$ 
the following two conditions are equivalent:
\begin{enumerate}
\item[(i)] $T$ is a finite multiplier sequence;
\item[(ii)] All different from $0$ roots of the polynomial 
$P_{T}(x)=\sum_{j=0}^kC_k^j 
\gamma _{j}x^j$
are  of the same sign.
\end{enumerate}
\end{theorem}
\vspace{1mm}

We get  by Theorem~\ref{th:MS} a linear diffeomorphism of 
$\mathcal M_n$ and $\overline {Hyp}_n^+\cup \overline{Hyp}_n^-$ 
where $\overline X$ means the closure of $X$.  
In this note we study the relation between the root 
multiplicities of $P$, $Q$ and $P\ast Q$.  
\vspace{1mm}

\begin{e-proposition}\label{th:mult}Given two (complex) polynomials $P$ 
and $Q$ of degree $n$  such  that $x_P$, $x_Q$ are roots respectively 
of $P$, $Q$ of multiplicity 
$m_P$, $m_Q$ with $\mu ^*:=m_P+m_Q-n\ge 0$, 
one has that $-x_Px_Q$ is a root of $P\ast Q$ of multiplicity 
$\mu ^*$. (If  $\mu ^*=0$, 
then $-x_Px_Q$ is not a root of  $P\ast Q$.)
\end{e-proposition} 
\vspace{1mm}

\begin{remark}
If $m_P>0$, $m_Q>0$ and $\mu ^*<0$, then $-x_Px_Q$  might or might not be a 
root of $P\ast Q$. Example:  $((x-1)(x-2)(x-3))\ast ((x-1)(x-4)(x-d))$ 
has $-1$ as a root if and only if $d=17/23$.
\end{remark}
\vspace{1mm}

\begin{e-proposition}\label{prop:comp}  For any 
%non-negative triple $u+v+w=n$, 
$P\in H_{u,v,w}$ and any $Q\in Hyp_n^-$ one has   
$P\ast Q\in H_{u,v,w}$. In particular, $Hyp_n^-$ is a semigroup w.r.t. the  
Schur-Szeg\"o  composition.   

\end{e-proposition}
\vspace{1mm}

The roots of $P$, $Q$ and $P\ast Q$ involved in  Proposition~\ref{th:mult} 
(i.e. 
those of the form $-x_Px_Q$, the sum of the multiplicities of $x_P$ and $x_Q$ 
being $>n$) are called {\em $A$-roots}, the remaining roots  
of $P$, $Q$, $P\ast Q$ are called $B$-roots.  With one exception -- if  
$0$ is a root of $P$, then it is considered as A-root of $P\ast Q$.   
Associate to $P\in Hyp_n$ its {\em multiplicity vector $MV_P$}
(the ordered partition of $n$ defined by the multiplicities of the  roots of 
$P$ in the 
increasing order).  For a root $\alpha$ of $P\in Hyp_n$ 
denote by 
$[\alpha ]_-$ (resp. $[\alpha ]_+$) the total number of roots of $P$ 
to the left 
(resp. to the right) of $\alpha$ and by sign$(\alpha )$ the sign of $\alpha$. 
 
%The main result of this note strengthens Proposition~\ref{prop:comp}.  
\vspace{1mm}

\begin{theorem}\label{th:ordpart} 
 For any $P\in Hyp_n$ and $Q\in Hyp_n^-$ the multiplicity 
vector $MV_{P\ast Q}$ 
is uniquely determined by Proposition~\ref{prop:comp} 
and the following conditions: 
\begin{enumerate}
\item[(i)]  For any A-root $\alpha \neq 0$ of $P$ and any A-root $\beta$ 
of $Q$ 
one has $[-\alpha \beta ]_-=[\alpha ]_-+[\beta ]_{sign(\alpha )}$. 
%where $sign(\alpha )$ is the sign of $\alpha$. 
\item[(ii)] Every $B$-root of $P\ast Q$ is simple. 
\end{enumerate}
\end{theorem}
\vspace{1mm}

\begin{corollary}The Schur-Szeg\"o composition restricted to $Hyp_n^-$ 
induces a semigroup structure on the set of all ordered partitions of $n$. 
Examples:   
$(2,14,1)\ast (5,6,6)=(1,1,2,1,1,1,3,1,1,1,3,1)$,      
$(1,14,2)\ast (5,6,4,2)=(1,2,1,1,1,3,1,1,1,1,1,1,1,1).$ 

\end{corollary} 

%\vspace{1mm}
%\noindent 
%Example.   $(2,14,1)\ast (5,6,6)=(1,1,2,1,1,1,3,1,1,1,3,1);\quad $     
%$(1,14,2)\ast (5,6,4,2)=(1,2,1,1,1,3,1,1,1,1,1,1,1,1).$ 
\vspace{1mm}

%\begin{remark}
%The $n$-tuple of reals $X=(x_1\le x_2\le ... \le x_n)$ is said 
%to be smaller than    $Y=(y_1\le y_2\le ... \le y_n)$ 
%{\em in the spectral order} (notation $X\preceq Y$) if their sums 
%coincide and for any $1\le i \le n$ 
%one has $x_1+x_2+\cdots +x_i \ge y_1+y_2+\cdots +y_i$, 
%see  e.g. \cite{HLP}. It 
%is interesting to know whether 
%\begin{conjecture}%[comp. \cite{Ju}] 
%composition in $Hyp_n^-$ preserves the spectral order, i.e. 
%if for $P_1,P_2,Q\in Hyp_n^-$ and 
%$P_1\preceq P_2$ one has $P_1*Q\preceq P_2*Q$. And also for  
%\end{conjecture}
%different $u_1+v_1+w_1=n$ and $u_2+v_2+w_2=n$ what is the 
%minimal/maximal number of real zeros of $P\ast Q$ where 
%$P\in H_{u_1,v_1,w_1},\;Q\in H_{u_2,v_2,w_2}$ (questions   
%raised in discussions with J.~Borcea). 
%\end{remark}
  
%\begin{ack}The authors are grateful to Wenner-Grenn foundation 
%for the support 
%of the visit of V.K. to Stockholm and to  J.~Borcea for  discussions of  
%the notion of  spectral order. 
%\end{ack} 

\section {Proofs.}

\noindent 
{\em Proof of Proposition~\ref{th:mult}:} Let $x_P$ and $x_Q$ be the required 
roots of $P$ and $Q$. 
It suffices to consider the case  $x_P =x_Q =-x_Px_Q=-1$. Indeed, 
$x_P$, $x_Q$, $-x_P x_Q$ 
are roots of $P(x)$, $Q(x)$, $P\ast Q(x)$ of multiplicities $m_P$, $m_Q$, 
$\mu ^*$ if and
only if $1$, $1$, $-1$ are roots of $P(x_P x)$, $Q(x_Q x)$, 
$P\ast Q(x_P x_Q  x)$ of the same multiplicities. 
%Notice that $P(x)\ast Q(x_P x_Q x)=P(x_P x)\ast Q(x_Q x)$. 
%Our proposition follows immediately from the next 
%lemma whose proof is straightforward. 
Set $G(x)=x^nP(1/x)$. Hence, $1$ is an $m_P$-fold root of $G$. One has 
$G^{(\nu )}(1)=\frac{n!}{(n-\nu )!}\sum _{j=0}^{n-\nu }C_{n-\nu}^ja_j$, 
$Q^{(\nu )}(1)=\frac{n!}{(n-\nu )!}\sum _{j=0}^{n-\nu }C_{n-\nu}^jb_{j+\nu}$. 
Set $K_s:=\frac{s!}{n!}G^{(n-s)}(1)=\sum _{q=0}^s
C_s^qa_q$, $L_r:=\frac{r!}{n!}Q^{(n-r)}(1)=\sum _{q=0}^rC_r^qb_{q+n-r}$. 
Hence, 
$K_n=K_{n-1}=\ldots =K_{n-m_P+1}=0=L_n=L_{n-1}=\ldots =L_{n-m_Q+1}$, 
$K_{n-m_P}\neq 0\neq L_{n-m_Q}$. One has 
$\sum _{j=0}^n(-1)^jC_n^jK_jL_{n-j}=
\sum _{j=0}^n(-1)^jC_n^ja_jb_j=(P*Q)(-1)~(*)$. (Indeed, to prove the 
equality between two bilinear forms in $a_i$, $b_k$, it suffices to set 
$a_{i_0}=1$, $a_i=0$ for $i\neq i_0$, $i_0=0,1,\ldots ,n$. The middle 
part of $(*)$ then equals $(-1)^iC_n^ib_i$, one has $K_j=C_j^i$ for $j\geq i$, 
$K_j=0$ for $j<i$, the left side equals $\sum _{j=i}^n(-1)^jC_n^jC_j^i
\sum _{\nu =0}^{n-j}C_{n-j}^{\nu}b_{\nu +j}~(**)$ and one checks directly that 
the coefficient before $b_l$ in $(**)$ equals $(-1)^iC_n^i$ if $l=i$ and $0$ 
if $l\neq i$.) 
Hence, if $\mu ^*>0$, then $-1$ is a root of $P\ast Q$ -- each product in 
the left side of $(*)$ contains a zero factor. When $\mu ^*=0$, then all 
but one products contain such a factor, so $(P\ast Q)(-1)\neq 0$. 
To prove that 
$-1$ is a root of $P\ast Q$ of multiplicity $\mu ^*$ one has to show for 
$\lambda <\mu ^*$ that 
$\sum _{j=0}^{n-\lambda}(-1)^jC_{n-\lambda}^jK_j^{\lambda}
L_{n-\lambda -j}=
\sum _{j=0}^{n-\lambda}(-1)^jC_{n-\lambda }^ja_{j+\lambda}b_{j+\lambda}=
\frac{(n-\lambda )!}{n!}(P*Q)^{(\nu )}(-1)$ where 
$K_j^{\lambda}=\sum _{q=0}^jC_j^qa_{q+\lambda}=
\frac{j!}{n!}(x^{n-\lambda}P^{(\lambda )}
(1/x))^{(n-\lambda -j)}|_{x=1}$.~~~~~$\Box$\\ 
%$L_j^{\lambda}=\sum _{q=0}^jC_j^qb_{q+n-j+\lambda}$.
\vspace{1mm}

\noindent 
{\em Proof of Proposition~\ref{prop:comp}:} 
We prove it  in the case $v=0$, i.e for any  
$P\in Hyp _n$, $P(0)\neq 0$ and any $Q\in Hyp_n^-$. 
The general case follows by continuity. The statement is trivially 
true for any hyperbolic  $P(x)$ of degree $n$ and $Q(x)=(1+x)^n$ since 
$P\ast Q=P$.  Let $Q\in Hyp_n^-$. 
Connect $Q(x)$ to $(1+x)^n$ by some path $Q^t(x)$ within 
$Hyp_n^-$. (This is possible since $Hyp_n^-$  is contractible.) 
Notice that if $P(0)\neq 0$, then  
$P\ast Q^t(0)\neq 0$ for the whole family 
since the constant term of $P\ast Q$ is the 
product of the ones of $P$ and $Q$. Therefore the number of 
positive and negative roots of $P\ast Q$ is the same as 
for $P\ast (1+x)^n$. $\Box$
\vspace{1mm}

{\em Proof of Theorem~\ref{th:ordpart}(i):}  Instead of $(P*Q)(x)$ we consider 
$Z(x):=(P*Q)(-x)$ (to have the same ordering of the roots
on the line in all three polynomials).
Suppose that $P\in H_{s,l,n-s-l}$, i.e.
$P=(\prod _{j=1}^s(x+a_j))x^l(\prod _{j=s+l+1}^n(x-a_j))$, $a_j>0$. 
Fix $a_j$ for $j=1,\ldots ,s$ and deform them continuously into $0$ for 
$j=s+l+1,\ldots ,n$. The MV of the negative root sets of $Z$ does not
change. Therefore to find this MV it suffices to find it for $P$ replaced by 
$P_1:=(\prod _{j=1}^s(x+a_j))x^{n-s}$. In the same way, to find the
MV of the positive root sets of $Z$ it suffices to find it for $P$ replaced by 
$P_2:=x^{s+l}(\prod _{j=s+l+1}^n(x-a_j))$. When $P(x)$ is changed
to $P(-x)$, then $Z(x)$ changes to $Z(-x)$, (this explains the presence of 
sign$(\alpha )$ in (i)) and the description of the 
MV of the root sets of $Z$ can be done by considering only
polynomials of the form $P_2$.

Consider the case when $P, Q\in Hyp_n^+$ (hence, $P$ and 
$Q$ play the same role). The other three cases $P\in Hyp_n^{\pm}$, 
$Q\in Hyp_n^{\pm}$ can be treated by analogy using 
$P(-x)*Q(x)=P(x)*Q(-x)=(P*Q)(-x)$. If $P=(x-a)^n$, then $Z(x)=Q(ax)$, 
so assume that 
each polynomial $P$, $Q$ has two distinct roots, 
$0<a_1<a_2$ and $0<b_1<b_2$, of multiplicities $m_1,m_2$ and $n_1,n_2$. If
$n$ is even and $m_1=n_1=n/2$, then $Z$ has no A-roots. Recall that by 
Proposition \ref{th:mult} if $a_ib_j$ is an A-root, then its multiplicity is
$m_i+n_j-n$.

Assume that (one of) the biggest of the four multiplicities $m_1$, $m_2$, 
$n_1$, $n_2$ is among the last two. Suppose first that this is $n_1$. If
$n_1+m_1>n$, $n_1+m_2>n$, then the root set $R(Z)$ of $Z$ looks like this: 
$(a_1b_1,V,a_2b_1,Y)$, see Propositions~\ref{th:mult} and \ref{pr:hyper}. 
Set $\sharp (V)=v$, $\sharp (Y)=y$. When writing 
$b_1\rightarrow 0$ or $b_1\rightarrow b_2$ 
we mean that the roots $a_1$, $a_2$, $b_2$ are fixed.
When $b_1\rightarrow 0$, then in the limit $Z$ has $n_2$ non-zero roots which 
are all from $Y$, hence, $y\geq n_2$. When $b_2\rightarrow b_1$,
then in the limit $Z(x)=P(b_1x)$ has two roots, of multiplicities 
$m_1$ and $m_2$. Hence, $(n_1+m_1-n)+v\geq m_1$, i.e. $v\geq n_2$. But
$v+y=2n_2$, hence, $v=y=n_2$.

If $n_1+m_1>n\geq n_1+m_2$, then $R(Z)=(a_1b_1,V)$, $v=n_2+m_2$. If 
$n_1+m_1\leq n<n_1+m_2$, then $R(Z)=(U,a_2b_1,V)$, $u+v=n_2+m_1$. When
$b_1\rightarrow 0$, then $v\geq n_2$ because all $n_2$ non-zero (in the limit) 
roots are in $V$. When $b_2\rightarrow b_1$, then
$(n_1+m_2-n)+v\leq m_2$, i.e. $v\leq n_2$. Hence, $v=n_2$, $u=m_1$.

Let $n_2=\max (m_1,m_2,n_1,n_2)$. If $n_2+m_1>n$, 
$n_2+m_2>n$, then $R(Z)=(U,a_1b_2,V,a_2b_2)$. When $b_1\rightarrow 0$, 
this yields $u\geq n_1$, and $b_1\rightarrow b_2$ yields $v\geq n_1$. 
As $u+v=2n_1$, one has $u=v=n_1$. If $n_2+m_1>n\geq n_2+m_2$, then
$R(Z)=(U,a_1b_2,V)$. When $b_1\rightarrow 0$, this yields $u\geq n_1$, 
and $a_1\rightarrow 0$ implies $v\geq m_2$. As $u+v=m_2+n_1$, one has
$u=n_1$, $v=m_2$. If $n_2+m_1\leq n<n_2+m_2$, then $R(Z)=(U,a_2b_2)$, 
$u=m_1+n_1$. This proves (i) of Theorem~\ref{th:ordpart} 
for $P,Q\in Hyp _n^+$ having each $\leq 2$ distinct roots.

%Thus Theorem 4 and Corollary \ref{samenum} are proved (without the assertion 
%that all B-roots of $P*Q$ are simple) in the case when there are up
%to two distinct roots in $P$ and $Q$, all roots being positive.

Further we assume that $P$ has a single A-root $a$, of multiplicity $m$. 
To prove the theorem by induction on the number of distinct 
positive roots in $P$ and $Q$ it suffices to consider the 
result on the MV of the root sets of $Z$ when a multiple root of 
$P$ or $Q$ splits into two. If this is a B-root, such a splitting 
deforms continuously the B-roots in $Z$, its A-roots and their
multiplicities don't change, and the theorem holds.

If an A-root splits into two B-roots, then it is a root of $Q$. 
Suppose that this is $b_{j_d}$, and that one has
$b_{j_{d-1}}<b_{j_d}<b_{j_{d+1}}$, $b_{j_{\nu}}$ being A-roots. 
Denote the multiplicities of these three roots 
by $h_1$, $h_2$, $h_3$, and by $t_1$, $t_2$ the sums of the
multiplicities of the B-roots of $Q$ from  
$(b_{j_{d-1}},b_{j_d})$ 
and $(b_{j_d},b_{j_{d+1}})$. Before the splitting of
$b_{j_d}$ the polynomial $Z$ had three A-roots stemming from 
$b_{j_{d-1}}<b_{j_d}<b_{j_{d+1}}$, namely, 
$ab_{j_{d-1}}<ab_{j_d}<ab_{j_{d+1}}$, 
of multiplicities $m+h_i-n$, $i=1,2,3$, with sums of the multiplicities of the 
B-roots of $Z$ from the two intervals between them equal
to $t_1+n-m$, $t_2+n-m$. After the splitting there 
remain only the 
A-roots $ab_{j_{d-1}}<ab_{j_{d+1}}$, the A-root $ab_{j_d}$ splits
into B-roots of total multiplicity $m+h_2-n$. In $Q$ there remain the A-roots 
$b_{j_{d-1}}<b_{j_{d+1}}$ with total multiplicity of the B-roots
between them equal to $t_1+t_2+h_2$. Thus the sum of the multiplicities of the 
B-roots of $Z$ from $(ab_{j_{d-1}},ab_{j_{d+1}})$
after the splitting equals $t_1+t_2-2m+2n+m+h_2-n=t_1+t_2+h_2+n-m$. Hence, 
(i) of Theorem~\ref{th:ordpart} 
holds after the splitting. If the A-root $b_{j_d}$ is 
first or last, i.e. $d=1$ or $d=r$, then the proof is similar.

Suppose that an A-root (say, $c$ of $P$, of multiplicity $\mu$) is splitting 
into an A-root to the left and a B-root to the right, of
multiplicities $\xi$ and $\eta$. Then in $Z$ there is a splitting of an A-root 
$cf$ ($f$ is an A-root of $Q$) of multiplicity $\mu +\nu -n$ into
an A-root of multiplicity $\xi +\nu -n$ and one or several B-roots of total 
multiplicity $\eta$. Suppose that at least one of these B-roots goes
to the left. Shift to the left (after the splitting) all roots of $P$ 
simultaneously while keeping the ones of $Q$ fixed. When one has $c=0$, then
the number of positive roots (counted with the multiplicities) will be greater 
for $P$ than for $Z$ (this follows from $[cf]_+=[c]_++[f]_+$ before
the splitting). This is a contradiction with Proposition \ref{prop:comp}. 
Hence, all new B-roots of $Z$ go to the right after the splitting and one
checks directly that (i) of Theorem~\ref{th:ordpart} holds after the 
splitting. If the B-root of $P$ goes to the left, or if $c$ is a
root of $Q$, then the reasoning is similar.

If an A-root $c$ splits into two A-roots $c^1$ (left) and $c^2$ (right) 
(hence, $c$ is a root of $Q$), then the above reasoning 
shows that in $Z$ an A-root $cf$ splits into two A-roots $c^1f$ 
(left, $f$ is an A-root of $P$) and $c^2f$ (right) and one or several
B-roots between them. Indeed, one shows as above that 
all roots different from $c^1f$ (resp. $c^2f$) and stemming from
$cf$ must go right (resp. left). Hence, the B-roots of $Z$ resulting from the 
splitting are between $c^1f$ and $c^2f$. Denote by $n^0$, $n^1$,
$n^2$ and $m^0$ the multiplicities of $c$, $c^1$, $c^2$ and $f$ 
($n^0=n^1+n^2$). Hence, the multiplicities of $c^1f$, $c^2f$ and the
total multiplicity of the B-roots of $Z$ between them equal $n^1+m^0-n$, 
$n^2+m^0-n$ and $n-m^0$, i.e. (i) of Theorem~\ref{th:ordpart} 
holds after the splitting. 
\vspace{1mm}

{\em Proof of Theorem~\ref{th:ordpart}(ii):}  
We show that non-simplicity of a B-root contradicts 
$P\ast Q\in Hyp_n$ for any $P\in Hyp_n$, $Q\in Hyp_n^{\pm}$, see 
Proposition~\ref{pr:hyper}. We first settle the basic case when either 
$P$ or $Q$ has only simple zeros and then use a procedure which either 
decreases the multiplicity of some root of $P$ or leads to  
$P\ast Q\not\in Hyp_n$. The multiplicity of $0$ as a root of 
$P$ must decrease up  
to $0$, not to $1$. 
%We first consider the case when $P(0)\neq 0$. 

{1) Basic case. } Suppose that $b \neq 0$ is a B-root of $P\ast Q$ of 
multiplicity $\mu \geq 2$. If $P$ has distinct real non-zero roots, then such
are all polynomials from a small neighbourhood $\Delta$ of $P$ in 
$Pol_n^{\mathbb{R}}$.   If $\mu =2$, then adjusting  the 
constant term which is non-zero by assumption one can easily choose 
$T\in \Delta$ such that $T\ast Q$ have a complex conjugate pair 
of zeros close to $b$ -- a contradiction. If $\mu >2$, then one
can choose $T$ such that $(T\ast Q)'$ has a multiple root at $b$ and
$(T\ast Q)(b)\neq 0$. Hence, $T\ast Q\not\in Hyp_n$.

{2) General case.} Assume that 
%$P(x)$ has a multiple roots at $c$,  i.e. 
$P=(x-c)^lP_1(x),\;l\ge 2$, $P_1(c)\neq 0$, $c\neq 0$. 
Set $P_c(x):=P(x)/(x-c)$.  Consider the family of hyperbolic polynomials 
$P^{\delta}= P+\delta P_c,\;\delta \in \mathbb{R}~~(\dag)$. 
In this family the $l$-tuple  root $c$ splits  into the $(l-1)$-tuple root 
$c$ and an extra root $c-\delta$ which is simple unless it coincides with 
some other root of $P$. Set $U_c:=P_c\ast Q$. 
Further considerations split into $3$ subcases  below:
%(2.i) -- 2.iii)):  
%i) $(P_c*Q)(b)\neq 0$; ii) $(P_c*Q)(b)=(P_c*Q)'(b)=0$ and iii) 
%$(P_c*Q)(b)=0\neq (P_c*Q)'(b).$

{2.i)} If $U_c(b)\neq 0$, then we find a value of 
$\delta$ such that 
$P^{\delta}\ast Q\not\in Hyp_n$ -- a contradiction. Indeed,  
if $\mu$ is even, then choosing 
sign$(\delta)$ one obtains  that $P^{\delta}*Q$ has no real roots 
close to $b$. If $\mu$ is odd, choose its
sign so that the total multiplicity of the roots of $P*Q$ close to $b$ is 
$<\mu$ (when $U_c'(b)\neq 0$, then $P*Q$ can be made monotonous
close to $b$; when $U_c'(b)=0\neq U_c''(b)$, then choose
sgn$(\delta )$ so that there is a local minimum (maximum) of $P*Q$ close
to $b$ where $P*Q$ is positive (negative); if
$U_c'(b)=U_c''(b)=0$, then for $\delta \neq 0$, $b$ is a
degenerate critical point and a non-root of $P*Q$, hence, $P*Q\not\in Hyp_n$). 

{2.ii)} If  $U_c(b)=U_c'(b)=0$, then 
%in the family $(\dag)$ 
$P^{\delta}*Q$ still has a multiple B-root at $b$ 
and lower multiplicity of $c$. 
%We reduce the multiplicities of roots of $P$. 

{2.iii)} Suppose that  $U_c(b)=0\neq U_c'(b)$. Assume  first that this 
happens  for at least two distinct roots $c$, $d$ of $P$  of which $c$ is 
multiple or $c=0$. Set $P_{cd}:=P/((x-c)(x-d))$. Consider 
the $2$-parameter family of hyperbolic polynomials  
$P^{\delta , \varepsilon}(x)= P(x)+\delta P_c(x)+\varepsilon
P_d(x)+\delta \varepsilon P_{cd}(x)~~(\ddag)$. In this family   the root $c$ 
(and also $d$ when multiple) splits as in $(\dag)$. Observe that
$P_{cd}=(P_c-P_d)/(c-d)~~(+)$. Hence, $P_{cd}\ast Q(b)=0$. 
Set $\delta =-\varepsilon
(P_d\ast Q)'(b)/((P_c+\varepsilon P_{cd})\ast Q)'(b)$. 
For $\varepsilon \neq 0$ small enough one
has $((P_c+\varepsilon P_{cd})\ast Q)'(b)\neq 0$. 
With this choice of $\delta$ 
the polynomial $P^{\delta , \varepsilon}*Q$ still has a multiple root at 
$b$ and lower multiplicity of $c$.

{2.iv)} To finish the  argument notice that  
the only case to consider when  one 
cannot perform splittings of roots of $P$ is  when $P$ has a single 
multiple or zero root $c$ (of multiplicity $\nu$) with 
$U_c(b)=0\neq U_c'(b)$, 
and the remaining non-zero roots  $d_i$ of $P$ are all simple 
with $(P_{d_i}*Q)(b)=(P_{d_i}*Q)'(b)=0$. The same must be  true for $Q$; 
denote by $g$ the root of $Q$ of
multiplicity $\lambda >1$. But  then a suitable linear combination of $P$ 
and $P_{d_i}$ equals $(x-c)^n$
(see $(+)$ etc.). Hence, $(x-c)^n*Q=Q(cx)$ has a multiple root at $b$,
i.e. $b=cg$. As $b$ must be a B-root (and not an A-root) of $P*Q$, one
must have $\nu +\lambda \leq n$. If $\nu +\lambda =n$, then $b$ is a 
non-root of $P*Q$ by Proposition~\ref{th:mult} -- a contradiction. 
If $\nu +\lambda <n$, then a suitable linear combination of
$P$ and $P_{d_i}$ equals $Y:=(x-h)^{\nu}(x-c)^{n-\nu }$ for which
one has that $b$ is a multiple root of $Y*Q$ -- a contradiction with   
Proposition~\ref{th:mult} again. $\Box$

\end{document}